\documentclass{amsart}

\usepackage{amsmath,amsthm}
\usepackage{amsrefs,alltt}

\usepackage{verbatim}
\newenvironment{details}{\comment}{\endcomment}
\newcommand{\ShowDetails}
{%
\renewenvironment{details}
  {%
   \par\noindent%
   \makebox[0in][r]
   {\raisebox{-1.5ex}[0in][0in]{$\downarrow$}\kern.4em}\kern-.7em%
   \hrulefill\\
   }
   {%
   \par\noindent\makebox[0in][r]
   {\raisebox{.5ex}[0in][0in]{$\uparrow$}\kern.4em}\kern-.7em%
   \hrulefill\par\noindent%
   }
}

\newcommand{\bilin}{B}

\newcommand{\map}[5][t]{\begin{map_machine}{#1}{#2}{#3}%
   {#4}{#5}\end{map_machine}}

\newenvironment{map_machine}[5]
   {
    
    \begin{array}[#1]{rcl}
   #2 & {-}\!\!{-}\!\!{\longrightarrow} & #3\\ 
   #4 & \mapstochar\kern-.25em{-}\!\!{\longrightarrow} & #5
   }{\end{array}}

\renewcommand{\epsilon}{\varepsilon}

\newcommand{\forwards}{\mbox{``$\Rightarrow$''}}

\newcommand{\define}[1]{\textbf{#1}}

\newcommand{\one}{^{(1)}}
\newcommand{\two}{^{(2)}}
\newcommand{\three}{^{(3)}}
\newcommand{\restrict}[1]{\big|_{\scriptstyle #1}^{}}

\renewcommand{\k}{k}

\newcommand{\R}{\mathcal{R}}

\DeclareMathOperator{\orthog}{O}
\renewcommand{\O}{\orthog}

\DeclareMathOperator{\GL}{GL}
\DeclareMathOperator{\Cl}{Cl}
\DeclareMathOperator{\symplectic}{Sp}
\def\Sp{\symplectic}
\DeclareMathOperator{\SO}{SO}

\newtheorem{theorem}{Theorem}[section]
\newtheorem{lemma}[theorem]{Lemma}
\newtheorem{cor}[theorem]{Corollary}
\newtheorem{prop}[theorem]{Proposition}

\theoremstyle{definition}
\newtheorem{definition}[theorem]{Definition}
\newtheorem{example}[theorem]{Example}
\newtheorem{remarks}[theorem]{Remarks}

\title[Unipotent classes and subgroups]
{Unipotent classes in the classical groups parameterized by subgroups}

\author{W. Ethan Duckworth}
\address{Loyola College, Baltimore MD\\
\texttt{educkworth@loyola.edu}}

\begin{document}
  
\maketitle

\begin{abstract}
  This paper describes how to use subgroups to parameterize unipotent
  classes in the classical algebraic group in characteristic $2$.
  These results can be viewed as an extension of the Bala-Carter
  Theorem, and give a convenient way to compare unipotent classes in a
  group $G$ with unipotent classes of a subgroup $X$ where $G$ is
  exceptional and $X$ is a Levi subgroup of classical type. \\
  \bigskip
  \textbf{AMS subject: } 14L35, 20G15
  \today
\end{abstract}

\section{Introduction and statement of results}
The unipotent classes in a classical group are often described using
Jordan blocks, but there are problems with this approach.  For
instance, the group $\SO_{2n}$, in all characteristics, and the groups
$\O_n$ and $\Sp_n$ in characteristic two, have distinct unipotent
classes with the same Jordan blocks sizes.  Furthermore, Jordan blocks
are defined using the natural module and are, therefore, not intrinsic
to the group.  For example, suppose that $G$ is an exceptional
algebraic group, and $X\le G$ is a Levi subgroup of classical type.
If we describe the unipotent classes in $X$ using Jordan blocks, then
this description does not make it easy to translate the list of
$X$-classes into unipotent classes in $G$.

In good characteristic the Bala-Carter Theorem avoids the problems
just described by describing the unipotent classes in all simple
algebraic groups using (pairs of) subgroups.  However, Bala-Carter
does not hold in bad characteristic.

The main goal of this paper is to extend the results of Bala-Carter
for groups of type $B_n$, $C_n$ and $D_n$ in characteristic $2$ and to
use this to translate the list of $X$-classes into unipotent classes
in $G$.  

For the rest of the paper let $\k$ be an algebraically closed field of
characteristic $p\ge 0$.  All the groups discussed here will be linear
algebraic groups over $\k$.  Each of the groups $\O_n$, $\SO_n$, and
$\Sp_n$ has a natural module of dimension $n$ which possesses a
bilinear form and a quadratic form.  A subspace is
\define{nonsingular} if it has trivial radical with respect to the
bilinear form, and it is \define{nondegenerate} if $0$ is the only
element of the radical which maps to $0$ under the quadratic form.  A
subspace is \define{totally singular} if the restrictions of the
quadratic and bilinear forms to the subspace are both identically
zero.  We use the following standard conventions to distinguish
between certain isomorphic but non-conjugate subgroups of $\SO_{2n}$.
A subgroup of $\SO_{2n}$ denoted by $\GL_{m}$ induces the full general
linear group on a pair of disjoint totally singular $m$-spaces which are in
duality via the bilinear form.  A subgroup of $\SO_{2n}$ denoted by
$\SO_{m}$ induces the special orthogonal group on a non-degenerate
$m$-space.

A subgroup of $\GL_n$, $\O_n$, or $\Sp_n$, denoted by $\Cl_m$, equals
$\GL_m$, $\O_m$, or $\Sp_m$, respectively.  We denote by $\Cl_m^\circ$
the identity component of $\Cl_m$.  Recall that $\GL_n$ and $\Sp_n$
are connected.  For the orthogonal groups we have $\O_m^\circ =
\SO_m$.  For $m\ge 2$, the group $\O_m$ equals $\O_m^\circ$ if and
only if $m$ is odd and $p = 2$.  The notation $\O_1$ means the trivial
group; we view it as acting on a $1$-dimensional vector space.

Let $G$ be one of $\GL_n$, $\SO_n$, or $\Sp_n$.  Let $\R(G)$
(``$\R$'' is chosen to stand for ``regular'') be the set of
closed subgroups $X \le \Cl_n$ where $X$ is a direct product of the
following form:
\begin{enumerate}
\renewcommand{\labelenumi}{(\roman{enumi})}
\item If $G= \GL_n$ then  $X = \GL_{n_1}\!\cdots \GL_{n_s}$ with $n_1
  + \dots + n_s = n$,
  
\item If $G = \SO_n$ and $p\ne 2$ or $G= \Sp_n$ then $X =
  \GL_{n_1}\cdots \GL_{n_s} \Cl_{m_1}^\circ \cdots \Cl_{m_r}^\circ$
  with $2n_1 + \cdots + 2n_s + m_1 + \cdots + m_r = n$,

\item If $G = \SO_n$ and $p=2$ then $X = \GL_{n_1}\cdots \GL_{n_s}
  \Cl_{m_1} \cdots \Cl_{m_r}$ with $2n_1 + \cdots + 2n_s + m_1 +
  \cdots + m_r = n$; if $n$ is even then $r$ even.

\end{enumerate}
\begin{details}
In each case there are two things that we want to be able to verify.
(1) that a regular unipotent element in $X$ is a unipotent element in
$G$ and (2) that we can get all the unipotent elements in $G$ in this
manner.

For (1) we note in case (i) or (ii) that $X$ is a subgroup of $G$.  In
case (iii) $X$ is not, but then we only need to verify that the
unipotent element have an even number of Jordan blocks, which we
ensure with the definition of $\R(G)$ below.  

For (2) we consider case-by-case.

  For $G=\SO_n$ with $p\ne 2$, we need to be able to get all
  partitions such that each even part has even multiplicity.  All
  parts with even multiplicity will come from regular elements in the
  $\GL_{n_i}$.  Furthermore, for odd parts with odd multiplicity we'll
  put all but one of the blocks into the $\GL_{n_i}$.  This leaves
  some number of distinct odd parts, call these $(m_1,\dots,m_r)$.
  Then the regular element in $\Cl_{m_i}^\circ$ gives $m_i$.  
  
  For $G=\Sp_n$ we do the same things, pull as many parts as necessary
  into the $\GL_{n_i}$.  What is left over is: $p\ne 2$, some distinct
  even parts, $p=2$ some even parts with multiplicity at most $2$.
  Call the left over parts $(m_1,\dots,m_r)$ etc.  
  
  For $G=\SO_n$ and $p=2$ we have the same parity conditions as
  $\Sp_n$.  Again pull as many parts as possible into the
  $\GL_{n_i}$.  What is left over is some even parts with multiplicity
  at most $2$.  We cannot get these from regular elements in
  $\SO_{m_i}$ with $m_i$ even, so we use $\O_{m_i}$.
\end{details}

If $X\in \R(G)$ we write $X =X_1 X_2 \cdots $ where the $X_i$ are the
factors of $X$ given in the definition of $\R(G)$.  The factors $X_i$
and the sequences $(n_1$, \ldots, $n_s)$ and $(m_1$, \ldots, $m_r)$
are uniquely determined (up to permutation) by $X$ and the
definitions.

A unipotent element in a connected reductive group $G$ is
\define{regular} if the dimension of its centralizer equals the rank
of $G$.  A connected reductive group has a single class of regular
unipotent elements \cite{steinberg_regular_elements}.  When $p=2$ and
$n$ is even we also consider the group $\O_n$ to have a unique regular
unipotent class in the non-identity component: an element in this
class also has dimension of its centralizer equal to the rank.  (See
Table~\ref{jordan_blocks_regular_unipotent_elements} for the Jordan
blocks of these classes; see \cite[I.4.8]{spaltenstein} for more
information about regular classes in non-connected reductive groups;
see \cite{hesselink1} for dimension of centralizer formulas).

\begin{theorem}
\label{theorem:map_with_regular_blocks}
Let $G$ be one of $\GL_n$, $\SO_n$, or $\Sp_n$.  Let $G$ act on
$\R(G)$ via conjugation.  The following map is surjective
\[
\Psi_1: \map {\{G\text{-classes in }\R(G)\}}
       {\{\text{unipotent }G\text{-classes}\}}
       {[2\jot]X = X_1 X_2 \cdots}
       {\text{conjugacy class of }u_1u_2\cdots}
\]
where each $u_i$ is a regular unipotent element in $X_i$, in the
non-identity component when possible.
\begin{details}
Thinking in terms of Jordan blocks, this is saying that all the
unipotent classes can be obtained by partitions satisfying the
restrictions imposed so far.
\end{details}

 Define a right inverse $\Phi_1$ of $\Psi_1$ as follows.  Given a
  unipotent $G$-class $C$, let $\Phi_1(C)$ equal the unique $G$-class
  in $\Psi_1^{-1}(C)$ which has a maximal number of factors of the form
  $\GL_{n_i}$.  The image of $\Phi_1$ equals all the $G$-classes in
  $\R(G)$ which satisfy the following conditions:
\begin{enumerate}
\renewcommand{\labelenumi}{(\roman{enumi})}
\item if $G=\GL_n$, then all of $\R(G)$ is in the image,
\begin{details}
In terms of Jordan blocks this says that the unipotent classes are
in bijection with all partitions.
\end{details}
  
\item if $G= \SO_n$ and $p\ne 2$, then the sequence $(m_1,\dots, m_r)$
  has distinct, odd parts,
\begin{details}
In terms of Jordan blocks this says that the unipotent classes are in
bijection with pairs of partitions $\alpha$, $\beta$ where $\alpha
\oplus \alpha \oplus \beta$ is a partition of $n$ and $\beta$ has
distinct, odd parts.  Compare this with \cite{carter}*{p394ff}
\end{details}

\item if $G= \Sp_n$ and $p\ne 2$, then the sequence $(m_1,\dots,m_r)$
  has distinct even parts,
  
\item if $p=2$ and $G$ equals $\SO_n$ or $\Sp_n$, then at most one part of
  the sequence $(m_1,\dots, m_r)$ equals $1$, and the rest of the
  parts are even with multiplicity at most $2$.
\begin{details}
In terms of Jordan blocks again, this says we have a pair $\alpha$,
$\beta$, with $\alpha \oplus \alpha \oplus \beta$ a partition of $n$,
$\beta$ has at most one part equal to $1$ and the rest of the parts
are even with multiplicity at most $2$.  Again compare with
\cite{carter}*{p394}.  
\end{details}
\end{enumerate}
\end{theorem}

\begin{remarks}
\label{remarks:after_theorem_maping_subgroups_to_regular_blocks}
(a) If $n$ is even, $G =\SO_n$, and $p=2$, then the groups in
$\R(G)$ are not always subgroups of $G$.  For the statement
of Theorem~\ref{theorem:map_with_regular_blocks}, this is unavoidable.
For example, when $p=2$, there is a distinguished unipotent class in
$\SO_{16}$ with Jordan block sizes given by $6$, $4$, $4$, $2$ and
this class cannot be represented by regular elements in a subgroup
with factors of the form $\GL_{n_i}$ and $\SO_{m_i}$ (c.f.
Remark~\ref{remarks:surjective_map_with_distinguished_parabolics} and
Example~\ref{example:classes_in_D_8}).
\begin{details}
  There are two classes in $\SO_{16}$ with partition given by
  $(6,4,4,2)$.  Using Spaltenstein's notation one class has
  $\epsilon(4) = 0$ and the other class has $\epsilon(4)=1$.  If the
  part $4$ comes from a regular element in $\GL_4$ then
  $\epsilon(4)=0$.  Therefore to get the class with $\epsilon(4)=1$
  the part $4$ has to come from an orthogonal group.  There is no way
  to get $(6,4,4,2)$ from regular elements of groups of the form
  $\SO_{m_i}$.
  
  This class is also not equal to the Richardson class of any
  parabolic subgroup.  This shows the need to use more than one
  parabolic subgroup in
  Theorem~\ref{theorem:surjective_map_with_distinguished_parabolics}
  below.  There are other classes which cannot be represented using
  two distinguished parabolic subgroups.  I do not know if one
  enlarged the family of parabolics if all classes could be obtained
  using only two parabolics.  
\end{details}

(b) Remark~\ref{remark:map_psi_1_is_properly_defined} shows that when
$n$ is even, $G=\SO_n$ and $p=2$, we have that $\Psi_1(X)$ is in $\SO_n$ and
not just in $\O_n$.
\end{remarks}

The next result is similar to
Theorem~\ref{theorem:map_with_regular_blocks}, except that we use
Richardson classes of all distinguished parabolic subgroups, instead
of using only regular classes.  Carter \cite{carter} provides a list
of distinguished parabolic subgroups for good characteristics, and we
use his list even in bad characteristics (note that his second formula
for $D_n$ has a slight mistake).
\begin{details}
  In Carter's notation from original paper \emph{Unipotent Elements in
    simple algebraic groups I} he had $m+n_1 + \dots + n_k = \ell -1$
  where $\ell$ is the rank.  Thus, $n_k$ is the horizontal length of
  the end nodes; or, if you like, the number of end nodes marked with
  a 0.  This number is correctly set equal to $\frac 12(n_{k-1}+1)$ or
  $\frac 12 n_{k-1}$.

However in the book Carter writes $m+ n_1 + \dots + n_k = \ell$.  Thus
$n_k$ is now used to denote the number of end nodes \emph{including}
the last $2$.  In this case he should set $n_k$ equal to $1 + \frac
12(n_{k-1}+1)$ or $1 + \frac 12 n_{k-1}$.  
\end{details}

Let $D(G)$ (``$D$'' standing for ``distinguished'') be the set of
closed connected subgroups $P$ of $G$ such that $P\le X$ for some
$X\in \R(G)$ (using notation as in the definition of $\R(G)$\,) with
the following changes: (1) if $G=\SO_n$ with $n$ even and $p=2$ we do
not require that $r$ be even; (2) we require $r\le 3$ and that if we
factor $X = X_1 X_2 \cdots$ then $P$ equals a direct product $P =
P_1P_2\cdots$ where for each $i$ we have that $P_i$ is a distinguished
parabolic subgroup of $X_i^\circ$.  The factors $P_i$ are uniquely
determined (up to permutation) by $P$.

In a connected reductive group, each parabolic subgroup has a unique
dense orbit in its unipotent radical \cite{richardson} which we call
the \define{Richardson class}.



\begin{theorem}  
\label{theorem:surjective_map_with_distinguished_parabolics}
Let $G$ be one of $\GL_n$, $\SO_n$, or $\Sp_n$.  Let $G$ act on $D
(G)$ via conjugation.  The following map is surjective,
\[
\Psi_2: \map {\{G\text{-classes in }D(G)\}}
     {\{\text{unipotent }G\text{-classes}\}}
     {[2\jot]P = P_1 P_2\cdots }
     {\text{conjugacy class of }u=u_1 u_2 \cdots}
\]
where $P$ is described above and each $u_i$ represents the
Richardson class of $P_i$.
\end{theorem}

\begin{remarks}
\label{remarks:surjective_map_with_distinguished_parabolics}
(a) The restriction of $\Psi_2$ to the subset of $D(G)$ consisting of
those $P$ with $r\le 1$ is injective and may be identified with the
map in the Bala--Carter Theorem from (pairs of) subgroups to unipotent
classes (see \cite{carter} and \cite{duckworth}).

(b) The distinguished unipotent class in $\SO_{16}$ mentioned in
Remark~\ref{remarks:after_theorem_maping_subgroups_to_regular_blocks}
does not equal the Richardson class of any parabolic subgroup.
\begin{details}
  See \cite[II.11.7]{spaltenstein} or reason as follows.  Let $\lambda
  = (6,4^2,2)$.  Then $\lambda^* = (4^2,3^2,1^2)$.  If you try to
  solve $\psi(\Lambda) = (4^2,3^2,1^2)$, where $\Lambda$ is the
  partition of the rank of $G$ arising from a Levi factor, then you
  find that $\Lambda$ must contain some $3$'s, but then
  $\psi(\Lambda)$ must contain $2$'s, a contradiction.
\end{details}
This indicates the need  in
Theorem~\ref{theorem:surjective_map_with_distinguished_parabolics} 
for more than one parabolic factor in $P$.  
\begin{details}
  As is shown below it is also not possible to use only two parabolic
  factors if we restrict ourselves distinguished parabolics (using
  only those parabolics on Carter's list).  However, there are two
  other possibilities: (1) enlarge the list, for example for $p=2$ and
  $C_n$ why not also use the parabolics which arise from $B_n$?  I
  haven't thought of an analogous way of enlarging the list for
  $D_n$.  (2) Use Spaltenstein induction instead of Richardson
  classes.  Then for classical groups one can obtain all the unipotent
  classes.  See \cite[II.7.10]{spaltenstein} ``In the classical
  groups, there are no rigid distinguished unipotent classes, except
  if $G^0$ is a torus.''
\end{details}

(c) When we define below the right inverse $\Phi_2$ of the map
$\Psi_2$ we refer to partitions consisting of Jordan blocks sizes.  In
principle, one can avoid mentioning Jordan blocks and still describe a
subset of the domain of $\Psi_2$, upon which $\Psi_2$ is injective,
thus implicitly describing an inverse of $\Psi_2$.  But the
description so obtained seems less natural than the presentation of
$\Phi_2$ given below.
\end{remarks}

A partition of $n$ is a sequence of natural numbers which add to $n$.
We write a partition $\beta$ as $(\beta_1$, $\beta_2$, $\dots)$ and
assume that $\beta_i \ge \beta_{i+1}$ for all $i$ unless otherwise
indicated.  We call each $\beta_i$ a part of $\beta$ (however we
sometimes have to keep track of the index $i$ in addition to the value
of $\beta_i$, see below).  Let $\alpha$ and $\beta$ be partitions of
$m$ and $n$ respectively.  We define $\alpha\oplus \beta$ to be a
partition of $m+n$ obtained by taking the union of the parts, counting
multiplicity, of $\alpha$ and $\beta$.

\begin{definition}
\label{remark_defining_decomposition}
Let $G$ equal $\SO_n$ or $\Sp_n$ and let $\beta$ be a partition.  We
will define a decomposition $\beta = \beta\one \oplus \beta\two \oplus
\beta\three$.  If $p\ne 2$, we set $\beta\two = \beta\three =0$.  If
$p=2$ and $G=\Sp_n$ then we set $\beta\three =0$ and we define
$\beta\two$ by the requirement that all parts of $\beta\two$ be
distinct and that a part of $\beta$ equal a part of $\beta\two$ if and
only if the part has multiplicity greater than $1$ in $\beta$.
\begin{details}
  In other words, if $\beta$ has only parts with multiplicity $\le 2$,
  then $\beta\two$ gets exactly one copy of each part which has
  multiplicity $2$.
\end{details}

If $p=2$ and $G = \SO_n$ we first define a map $f$ from parts of
$\beta$ to $0$ or $1$.  We allow the abusive notation that $\beta_i =
\beta_{i+1}$ but that $f(\beta_i) \ne (\beta_{i+1})$; in these cases
the subscript in $\beta_i$ is implicitly part of the definition of $f$.

Define $f(\beta_1) = 1$.  Let $j$ be given such that $f(\beta_1)$,
\dots, $f(\beta_{j-1})$ have been defined, let $\beta_k$ be the last
of these parts which maps to $1$, let $\ell$ and $i$ be the number of
parts in $\beta_1$, \ldots, $\beta_{j-1}$ which map to $1$ and to $0$
respectively.  Define $f(\beta_j)$ as follows (where we allow
$\beta_{j+1}$, $\beta_{j+2}$, etc. to equal $0$):
$$
f(\beta_j) = 
\left\{
\begin{array}{lp{2in}@{\qquad}c}
0 & if $\ell$ is even and $\beta_k - \beta_j \le 2$ & (1)\\
0 & if $\ell$ is even, $i$ is odd, $\beta_{j+1}\in \{0,1\}$ & (2) \\
0 & if $\ell$ is even, $i$ is odd, $\beta_{j+1}-\beta_{j+3}\le 2$,  
$\beta_{j+3} \ne 0$, $\beta_j - \beta_{j+3} \ge 3$ &
\raisebox{-3ex}{(3)} \\ 
1 & in all other cases. 
\end{array}
\right.
$$
(The result $f(\beta_j)=1$ is meant to be the generic case, with
conditions (1), (2) and (3) viewed as exceptions; see
Example~\ref{example:decomposing_with_f}.)

Finally, when $p=2$ and $G=\SO_n$ we apply $f$ to $\beta$, let
$\beta\one$ and $\delta$ equal the pre-image of $1$ and $0$
respectively.  Apply $f$ to $\delta$, let $\beta\two$ and
$\beta\three$ be the pre-image of $1$ and $0$ respectively.

\end{definition}

In the following theorem the notation $|\beta^{(i)}|$ denotes the sum of
the parts of $\beta^{(i)}$.

\begin{theorem}
\label{theorem:inverse_of_map_with_distinguished_parabolics}
If $G=\GL_n$ then $\Psi_2$ is bijective.  Otherwise, we define a right
inverse $\Phi_2$ of $\Psi_2$ as follows.  Let $u\in G$ be unipotent,
let $L = \GL_{n_1}\cdots \GL_{n_s}\Cl_m^\circ$ be a minimal Levi
subgroup containing $u$ and factor $u = u_1 \cdots u_s u_0$ with $u_i
\in \GL_{n_i}$ for $1\le i \le s$ and $u_0\in \Cl_m^\circ$.  Let
$\beta$ be the Jordan blocks of $u_0$ in the natural module for
$\Cl_m^\circ$ and write $\beta = \beta\one \oplus \beta\two \oplus
\beta\three$ as in
Definition~\textup{\ref{remark_defining_decomposition}}.  Then $u$ is
contained in a subgroup of $L$ of the form $X = X_1\cdots X_{s+3}$
with $X_i = \GL_{n_i}$ for $1\le i \le s$ and $X_{s+i} =
\Cl_{|\beta^{(i)}|}^\circ$ for $1 \le i \le 3$.  We factor $u_0$
further such that $u = u_1 \cdots u_{s+3}$ where $u_i\in X_i$ for each
$i$.  For $1\le i \le s+3$ there exists a distinguished parabolic
subgroup $P_i$ of $X_i$, unique up to conjugacy, whose Richardson
class is represented by $u_i$.  We define $\Phi_2(u)$ to equal the
$G$-orbit of $P = P_1\cdots P_{s+3}$.
\end{theorem}

When $p\ne 2$, Theorems
\ref{theorem:surjective_map_with_distinguished_parabolics} and
\ref{theorem:inverse_of_map_with_distinguished_parabolics} are
equivalent to the Bala--Carter Theorem for the classical groups.  By
an \define{extra class} we mean one which is not parameterized by
Bala--Carter.  In the notation of
Theorem~\ref{theorem:inverse_of_map_with_distinguished_parabolics} a
class is extra if and only if $\beta \ne \beta\one$.
  

\begin{cor}
\label{corollary:split_unipotent_classes_and_Levis}
Two unipotent classes in $\SO_{n}$ are conjugate under $\O_{n}$, but
not under $\SO_{n}$, if and only if these classes correspond under the
map of Theorem~\textup{\ref{theorem:map_with_regular_blocks}} to a pair of
Levi subgroups in $\R(G)$ (or, under the map in
Theorem~\textup{\ref{theorem:inverse_of_map_with_distinguished_parabolics}},
to a pair of Borel subgroups in these Levi factors) which are also
conjugate under $\O_{n}$ but not $\SO_{n}$.
\end{cor}

We note that
Corollary~\ref{corollary:split_unipotent_classes_and_Levis} is well
known for $p\ne2$ as it follows from the usual Bala-Carter Theorem.

In Section~\ref{remark:jordan_blocks_of_rich_classes} we give an
explicit, combinatorial formula for the Jordan blocks of $\Psi_1(X)$
or $\Psi_2(P)$.  This formula is used to determine the parabolics $P_i$
in Theorem~\ref{theorem:inverse_of_map_with_distinguished_parabolics}.
Proposition~\ref{decomposing_distinguished_classes_into_Richardson_classes}
establishes certain canonical properties possessed by the
decomposition $\beta^{(1)} \oplus \beta^{(2)} \oplus \beta^{(3)}$.

\section{Jordan block parameterization of unipotent classes}
\label{section_on_preliminaries}
In this section we recall one method of parameterizing unipotent
classes in $G$, following \cite{spaltenstein}*{I.2.5ff} (though we
extend the method there to include the odd-dimensional orthogonal case
in characteristic $2$).

For the moment we fix $G$ equal to one of $\O_n$ or $\Sp_n$ defined
over an algebraically closed field $\k$ of characteristic $p\ge 0$ and
we denote by $\bilin$ the associated bilinear form.  Calculations
involving $u-1$ are made by viewing the natural module for $G$ as a
$\k[u]$-module.

Let $\delta$ equal $1$ if $G = \Sp_n$ or if $G = \O_n$ and $p=2$.  Let
$\delta$ equal $-1$ otherwise.
  
Let $u\in G$ be unipotent and let $\lambda$ equal the Jordan block
sizes of $u$.  For each part $x$ of $\lambda$ let $\epsilon_u(x)$ be
defined as follows: if $x$ is odd let $\epsilon_u(x) = -\delta$, if
$x$ is even then $\epsilon_u(x) \in \{0,\delta\}$ with $\epsilon_u(x)
= 0$ if and only if $p=2$ and $\bilin \big ( (u-1)^{x -1}v,\ v\big) =
0$ for all $v \in \ker (u-1)^x$.  Usually, if the element $u$ has been
fixed or is irrelevant, we will write $\epsilon$ instead of
$\epsilon_u$.  This gives a map $\Upsilon$ from unipotent classes in
$G$ to pairs $(\lambda,\ \epsilon)$.  The map $\Upsilon$ is injective.
  
Let $\lambda$ be a partition of $n$ and $\epsilon$ a map from the set
$\{ \lambda_i \mid i\ge 1\}$, to the set $\{-1,0,1\}$.  Then
$(\lambda,\ \epsilon)$ is in the image of $\Upsilon$ if and only if
two of the following are satisfied:
\begin{enumerate}
\renewcommand{\labelenumi}{(\roman{enumi})}
\item $G= \O_n$, $p\ne 2$, every even part of $\lambda$ has even
  multiplicity,
  
\item $G=\O_n$, $p= 2$ or $G=\Sp_n$, every odd part of $\lambda$
  strictly greater than $1$ has even multiplicity,
  
\item $p\ne 2$, $\epsilon(x) = -\delta$ if $x$ is odd and $\epsilon
  (x) = \delta$ if $x$ is even,
  
\item $p=2$, $\epsilon(x) = -1$ if $x$ is odd, $\epsilon(x) = 1$ if
  $x$ is even with odd multiplicity and $\epsilon(x) \in \{0,1\}$ if
  $x$ is even with even multiplicity.
\end{enumerate}

This completes the parameterization of unipotent classes in $\O_n$ and
$\Sp_n$.  Now we relate the unipotent classes of $\O_n$ to the classes
in $\SO_n$.

If $p\ne2$ or if $n$ is odd, then every unipotent class in $\O_n$ is
contained in $\SO_n$.  
\begin{details}
For $p\ne2$ all the unipotent elements of $\O_n$ are contained in
$\SO_n$. 

If $p=2$ and $n$ is odd then $\SO_n = \O_n$ (see, for example,
\cite{carter}*{p181}).
\end{details}

If $p=2$ and $n$ is even, then a unipotent element in $\O_n$ is
contained in $\SO_n$ if and only if it has an even number of Jordan
blocks.

We apply the definition of $\epsilon$ to $\SO_n$ without change.  Let
$u\in \SO_{n}$ be unipotent, let $C$ be the $\O_{n}$-class of $u$, and
let $\lambda$ be the Jordan block sizes of $u$.  Then $C$ equals a single
$\SO_{n}$-class, unless each part $x$ of $\lambda$ is even and
$\epsilon(x) \ne 1$, in which case $C$ forms two $\SO_n$-classes.
\begin{details}
If $p=2$ and $n$ is odd then $\O_n =\SO_n$ so no splitting.

If $p\ne 2$ then the split occurs if and only if all the parts and all
the multiplicities are even (see \cite{spaltenstein}*{p19}) in which
case we have $\epsilon(x) = -1$.

If $p=2$ and $n$ is even then the split occurs if and only if all the
parts and all the multiplicities are even and $\epsilon(x)\ne 1$ for
all the parts (see \cite{spaltenstein}*{p20}).  
\end{details}

\begin{remarks}
\label{remark:map_psi_1_is_properly_defined}
We pause to clarify one aspect of the map $\Psi_1$ from
Theorem~\ref{theorem:map_with_regular_blocks}.  If $G = \SO_n$ with
$n$ even and $p=2$ then some elements of $\R(G)$ are
subgroups of $\O_n$ but not subgroups of $G$.  However the regular
unipotent class which we have chosen in such subgroups is contained in
$\SO_n$, so the map is properly defined.
\end{remarks}

\section{Jordan 
  blocks of Richardson classes of distinguished parabolics.}
\label{remark:jordan_blocks_of_rich_classes}
Let $u$ be a regular unipotent element in $G$.  In
Table~\ref{jordan_blocks_regular_unipotent_elements} we describe the
possible Jordan blocks of $u$.
\begin{table}
\caption{Jordan blocks of regular unipotent elements}
\label{jordan_blocks_regular_unipotent_elements}
$$
\begin{array}{cp{2in}}
\text{Jordan blocks} & Groups and conditions\\
\hline
n      & $\GL_n$, $\Sp_n$,                                       \\
       & $\SO_n$ with $n$ odd and $p\ne 2$,                       \\
       & $\O_n$ with $u$ in the non-identity component, $n$ even and
       $p=2$\\[1ex]
n-1, 1 & $\SO_n$ with $n$ odd and $p=2$ or $n$ even and $p\ne 2$ \\[1ex]
n-2, 2 & $\SO_n$ with $n$ even, $n\ge 4$, $p=2$ (and $\epsilon(2)=1$
if $n=4$)\\[1ex]
1, 1   & $\SO_2$ with $p=2$
\end{array}
$$
\end{table}

Let $G$ equal $\SO_n$ or $\Sp_n$, let $m$ be the rank of $G$, and let
$P$ be a parabolic subgroup of $G$.  Let $L$ be a Levi factor for $P$
and write $L=\GL_{n_1} \cdots \GL_{n_s} H_{m_0}$ where $H_{m_0}$ is
one of $\SO_{2m_0}$, $\SO_{2m_0+1}$, or $\Sp_{2m_0}$ (we allow $m_0=0$
and $H_{m_0} = 1$).  In this manner
$P$ determines a partition of $m$ given by $m = n_1 + \cdots + n_s +
m_0$.  We write this partition as $(1^{c(1)}$, $2^{c(2)}$, $\dots$,
$N^{c(N)}_{})\oplus m_0$ where $N= \max\{n_i \mid 1\le i \le s\}$ and
$c(x)$ is the multiplicity of the part $x$ in the sequence $n_1$,
\dots, $n_s$.  We assume now that $P$ is distinguished.  Then we have
$c(i) \ge 1$ if and only if $1\le i \le N$.  ??? ref ???
Let $\lambda$ be the partition of $n$ whose parts equal the Jordan
block sizes of the Richardson class of $P$.  Recall that the dual of
$\lambda$ is the partition $\lambda^*$ of $n$ where $\lambda_i^*$
equals the number of $j$ such that $\lambda_j \ge i$.  Recall also
that $\lambda = (\lambda^*)^*$.  In
Table~\ref{table:levi_partitions_distinguished} we describe $\lambda$
in terms of its dual.  In this table we have written parts of
$\lambda^*$ which may have multiplicity $0$; for example, when $G =
\SO_{2m+1}$ and $p\ne 2$ we can have $c(2m_0+1)=0$.
Thus Table~\ref{table:levi_partitions_distinguished} implicitly
defines a map from $G$-classes of distinguished parabolics to
partitions of $n$.  This map is injective and its image is described
in Table~\ref{table_jordan_blocks_dist_rich_classes} (see
\cite{duckworth}).

Table~\ref{table:levi_partitions_distinguished} allows one to
calculate the Jordan blocks of $\Psi_2(P)$ where $P$ and $\Psi_2$ are
as in
Theorem~\ref{theorem:inverse_of_map_with_distinguished_parabolics}.
Table~\ref{table_jordan_blocks_dist_rich_classes} will allow us
verify, in Section~\ref{section_main_proofs}, the assertion in
Theorem~\ref{theorem:inverse_of_map_with_distinguished_parabolics}
that there exists a distinguished parabolic subgroup $P_i$ of $X_i$,
unique up to conjugacy, whose Richardson class is represented by
$u_i$.
\begin{table}
\small
\caption{\small Jordan blocks of the Richardson class of a
distinguished parabolic.}%
\label{table:levi_partitions_distinguished}%
\[
\begin{array}{crl}
\hline
G &   & \qquad  \lambda  \\ 
\hline
\GL_m & m_0 =0 &  (1^{c(1)},\ 2^{c(2)},\ \dots)^*=(m)=(n)\\
[2ex]
\SO_{2m} 
   &  
     m_0 = 0,\   p\ne 2: & \left(1^{2c(1)-2},\ 2^{2c(2)+1}\right)^*\\  
   & m_0 = 0,\   p=2: & \left(1^{2c(1)-4},\ 2^{2c(2)+2}\right)^*\\
   &  
    m_0>0,\ p\ne2: 
         & \left(1^{2c(1)},\ \dots,\ (2m_0-1)^{2c(2m_0-1)},\
         (2m_0)^{c(2m_0)+1}\right)^*\\ 
   & m_0>0,\  p=2:
         & \Big(1^{2c(1)-2},\ 2^{2c(2)+2}, 3^{2c(3)-2}, 4^{2c(4)+2}, \ \dots,\\
         && \hfill (2m_0-1)^{2c(2m_0-1)-2},\ (2m_0)^{2c(2m_0)+2}\Big)^*\\ 
[2ex]
\SO_{2m+1} 
   & 
       p\ne 2: 
         & \left(1^{2c(1)},\ 2^{2c(2)},\ \dots,\ (2m_0)^{2c(2m_0)},\
         (2m_0+1)^{2c(2m_0+1)+1}\right)^*\\ 
   &   p=2:
         & \Big(1^{2c(1)-2},\ 2^{2c(2)+2},\ 3^{2c(3)-2},\ 4^{2c(4)+2}, \dots,\\
         && \hfill (2m_0)^{2c(2m_0){+}2},\ (2m_0{+}1)^{2c(2m_0{+}1)}\Big)^*\oplus (1)\\ 
[2ex]
\Sp_{2m} 
   & m_0=0: 
         & \left(1^{2c(1)},\ 2^{2c(2)},\ \dots,\ N^{2c(N)}\right)^*\\
\hline
\end{array}
\]  
\end{table}

\begin{table}
\small
\caption{\small Partitions which equal Jordan blocks of the Richardson
class of a distinguished parabolic}%
\label{table_jordan_blocks_dist_rich_classes}%
\[
\begin{array}{lp{4in}}
\hline
G & Set of partitions $\lambda$ (\,$\epsilon(x)= 1$ for all parts $x$
unless otherwise noted)\\ 
\hline
G=\GL_m & All partitions of $m$\\[1\jot]

G= \SO_{2m+1},\ p\ne 2 & Partitions of $2m+1$ consisting of distinct 
odd parts\\[1\jot]
G=\SO_{2m+1},\ p=2 & Partitions of $2m+1$  which have exactly
one part equal to $1$, $\epsilon(1){=}{-}1$, the rest of the parts
are even and of multiplicity at most $2$, and, if $i$ 
is even and $\lambda_{i+1}\ge1$ then $\lambda_i - \lambda_{i+1} \ge
3$\\[1\jot] 
G=\Sp_{2m} & Partitions of $2m$ consisting of distinct even 
parts\\[1\jot]
G= \SO_{2m},\ p\ne 2 & Partitions of $2m$ consisting of distinct odd 
parts\\[1\jot]
G=\SO_{2m},\ p=2 & Partitions of $2m$ which have an even number of
parts, each part is even and of multiplicity at most $2$,
and, if $i$ is even and $\lambda_{i+1}\ge1$ then $\lambda_i -
\lambda_{i+1} \ge 3$\\ 
\hline
\end{array}
\]
\end{table}

\section{Main proofs}
\label{section_main_proofs}
Let $G$ be connected and reductive and let $L$ be a Levi subgroup of
$G$ (we allow $L=G$).  A unipotent element $u\in L$ is
\define{distinguished} (in $L$) if $u$ is not contained in any proper
Levi subgroup of $L$.  If we omit mention of $L$ then we assume $L=G$.
Then $u$ is distinguished if and only if each maximal torus of
$C_G^{}(u)$ is contained in $Z(G)$ (see
Lemma~\ref{min_Levi_unique_up_to_conj}).  If $Z(G)=1$ this is
equivalent to requiring that $C_G^{}(u)$ have no nontrivial torus (this
is the usual definition) which is also equivalent to requiring that
$C_G^{}(u)^0$ be a unipotent group.


For many questions, the following lemma reduces the study of unipotent
classes in $G$ to the study of distinguished classes.

\begin{lemma}[\cite{carter}*{5.9.2, 5.9.3}]
\label{min_Levi_unique_up_to_conj}
\begin{enumerate}
\renewcommand{\labelenumi}{(\roman{enumi})}
\item Let $S$ be a torus.  Then $L=C_G^{}(S)$ is a Levi subgroup.
  
\item If $u$ is a unipotent element and $S$ a maximal torus of $C_G^{}
  (u)$ then $u$ is distinguished in $L = C_G^{} (S)$.  Furthermore,
  any Levi subgroup in which $u$ is distinguished is conjugate to $L$
  via an element of $C_G^{}(u)^\circ$.
\end{enumerate}
\end{lemma}

\begin{cor}
\label{bijection_reduction_to_dist_classes} Define a map from
$G$-classes of pairs $(L,C)$ consisting of a Levi subgroup, $L$, of
$G$ and a distinguished unipotent $L$-class, $C$, to unipotent
$G$-classes by extending $C$.  This map gives a bijection.
\end{cor}

In the following lemma the notation $\lambda(u_i)$ for $i\ge 1$ denotes
the Jordan block sizes of $u_i$ in the natural module for
$\GL_{n_i}$.  The notations $\lambda(u)$ and $\lambda(u_0)$ denote the
Jordan block size of $u$ and $u_0$ in the natural module for $G$ and
$\Cl_m^\circ$ respectively.  

\begin{lemma} 
\label{lemma_mapping_subgroups_to_classes}
Let $G$ equal $\SO_n$ or $\Sp_n$ and let $u \in G$ be unipotent.  The
following hold:
\begin{enumerate}
\renewcommand{\labelenumi}{(\roman{enumi})}
\item Let $\GL_{n_1} \cdots \GL_{n_s}\Cl_m^\circ$ be a Levi subgroup
  of $G$ containing $u$ and write $u = u_1 u_2\cdots u_s u_0$ with
  $u_i \in \GL_{n_i}$ for $i\ge1$ and $u_0\in \Cl_m^\circ$.  Then
  $\lambda(u) = \bigoplus_{i\ge 1} \lambda(u_i)^2 \oplus \lambda(u_0)$
  where $\lambda(u_i)^2$ means that each part of $\lambda(u_i)$ has
  been doubled in multiplicity.
  
\item If $p\ne 2$ then $u$ is distinguished if and only if each part
  of $\lambda(u)$ has multiplicity $1$ (whence each part is odd if $G
  = \SO_n$ and each part is even if $G = \Sp_n$).
  
\item If $p=2$ then $u$ is distinguished if and only if at most one
  part of $\lambda(u)$ equals $1$, and each remaining part $x$ has
  multiplicity at most $2$ and $\epsilon(x)=1$ (whence $x$ is even).

\item Let $V$ be the natural module for $G$, let $x$ be a Jordan block
  of $u$ and suppose that $u$ stabilizes a decomposition $V = V_1\perp
  V_2$.  Let $u_1 = u\restrict{V_1}$ and $u_2 = u\restrict{V_2}$.
  Then $\epsilon_u(x)=1$ if and only if $\epsilon_{u_1}(x)=1$ or
  $\epsilon_{u_2}(x)=1$.
  
\item Let $n_1$ be the rank of $G$, let $\GL_{n_1}$ be a Levi factor
  of $G$, and let $u$ be distinguished in $\GL_{n_1}$.  Suppose that
  $n_1$ is even and that $p=2$.  Then $\epsilon(n_1)=0$.
\end{enumerate}
\end{lemma}

\begin{proof}
  Parts (i), (ii) and (iii) are in \cite{spaltenstein}*{II.7.10}.
  
  Part (iv) (sketch).  The crucial case is where $x$ is even and
  $p=2$, which we now assume.
\begin{details}
If $x$ is odd then $\epsilon_u(x) = \epsilon_{u_i}(x) = -\delta$.

If $x$ is even but $p\ne 2$ then $\epsilon_u(x) = \epsilon_{u_i}(x) =
\delta$.  

If $x$ is even and $p=2$ but the multiplicity of $x$ is odd then the
multiplicity of $x$ in one of the $V_i$ is also odd, and then
$\epsilon_u(x) = \epsilon_{u_i}(x)=\delta$.  
\end{details}
Writing any $v\in V$ as $v= v_1 + v_2$ with $v_i\in V_i$, it is easy
to show that 
\begin{multline*}
\bilin\big((u-1)^{x-1}v,\ v \big)=0 \quad \forall v\in \ker(u-1)^x\\
\iff 
\bilin\big((u_i-1)^{x-1}v_i,\ v_i\big)=0 \quad \forall v_i\in \ker
(u_i-1)^x \text{ for }i=1,2.
\end{multline*}
The result now follows from the definition of $\epsilon$ (c.f.
Section~\ref{section_on_preliminaries}).

Part (iv).  From part (i) we know that the multiplicity of $n_1$ is
$2$.  Since $u$ is \emph{not} distinguished in $G$ we have by part
(iii) that $\epsilon(n_1)\ne 1$.
\end{proof}

The following result is essentially equivalent to
Corollary~\ref{corollary:split_unipotent_classes_and_Levis}.

\begin{cor}
Let $n$ be even, let $u_1$ and $u_2$ be unipotent elements in
$G=\SO_{n}$ and for $i=1,2$ let $L_i$ be a minimal Levi subgroup of
$G$ containing $G$.  Then $u_1$ and $u_2$ are conjugate under $\O_n$
but not $\SO_n$ if and only if $L_1$ and $L_2$ are conjugate under
$\O_n$ but not $\SO_n$. 
\end{cor}

\begin{proof}
  Recall that an $\O_n$-class of Levi subgroups splits into two
  $\SO_n$-classes if and only if the class is represented by
  $\GL_{n_1}\cdots \GL_{n_s}$, with each $n_i$ even.  (This can be
  shown by viewing each Levi subgroup as a stabilizer of subspaces and
  then using standard arguments about the geometry of classical groups
  and Witt's Theorem, or by more abstract arguments about conjugacy of
  subgroups and root systems of algebraic groups.)
\begin{details}
Let $V$ be the natural module for $G$.  We view each Levi factor as
being the stabilizer in $G$ of a direct sum decomposition
$$
V = (V_{1}\oplus V_{1}^*) \perp \dots \perp (V_{s}\oplus V_{s}^*)
\perp \widetilde V
$$
where each $V_{i}$ and $V_{i}^*$ are totally singular, each sum
$V_{i}\oplus V_{i}^*$ is nonsingular, and $\widetilde V$ is
nonsingular.  Using Witt's Theorem one can show that $L_1$ and $L_2$
are conjugate under $\O_n$ if and only if each stabilizes a direct sum
decomposition with the same list of dimensions for each $\dim V_i$,
$\dim V_i^*$ and $\widetilde V$.  Thus $L_1$ and $L_2$ are conjugate
under $\O_n$ if and only if we can write each as $\GL_{n_1}\cdots
\GL_{n_s}\SO_{2m}$ for the same $n_i$ and $m$.  Now one applies the
same reasoning using conjugation by $\SO_n$ instead of $\O_n$, with
the following change: two totally $n$-dimensional subspaces are
conjugate under $\SO_n$ if and only if their intersection has even
codimension in each of these subspaces.  

Note: Bala-Carter have a little discussion about the conjugacy of
subsystems.  They cite Dynkin as having classified all such
conjugacy classes.  
\end{details}
Recall from Section~\ref{section_on_preliminaries}, that an unipotent
$\O_n$-class splits into two $\SO_n$-classes if and only if each
Jordan block size $x$ is even and satisfies $\epsilon(x) \ne 1$.

Suppose that $L_1$ and $L_2$ are conjugate under $\O_n$ but not under
$\SO_n$.  Then $L_1\cong L_2 \cong \GL_{n_1}\cdots \GL_{n_s}$ with
each $n_i$ even.  Then $\lambda(u_1) = \lambda(u_2) =
(n_1^2,n_2^2,\dots,n_s^2)$ and, by
Lemma~\ref{lemma_mapping_subgroups_to_classes} (iv) and (v), we
have $\epsilon(n_i) \ne 1$ for each $i$.  Therefore $u_1$ and $u_2$
are conjugate under $\O_n$ but not $\SO_n$.  

Suppose that $u_1$ and $u_2$ are conjugate under $\O_n$ but not under
$\SO_n$.  By Corollary~\ref{min_Levi_unique_up_to_conj}(ii) we have
that $L_1$ and $L_2$ are conjugate under $\O_n$.
\begin{details}
Let $u_1 = gu_2$ for $g\in \O_n$.  Then $L_1$ and $gL_2$ are both
minimal Levi subgroups of $u_1$.  Then an element of $C_G^{} (u_1)$ takes
$gL_2$ to $L_1$.
\end{details}

Now we claim that $L_1$ and $L_2$ cannot have any factor of the form
$\SO_{2m}$.  Otherwise $u_1$ and $u_2$ could each be written as a
product with one factor distinguished in $\SO_{2m}$.  By
Lemma~\ref{lemma_mapping_subgroups_to_classes} and
Table~\ref{table_jordan_blocks_dist_rich_classes} this would give rise
to at least one part $x$ of $\lambda(u_1)$ and $\lambda(u_2)$ with
$\epsilon(x)=1$.  But then $u_1$ and $u_2$ would be conjugate under
$\SO_{2n}$.

Now $\lambda(u_1)$ and $\lambda(u_2)$ are equal (since $u_1$ and $u_2$
are conjugate under $\O_n$) with parts $(n_1^2,\dots, n_s^2)$ where
$L_1 \cong L_2 \cong \GL_{n_1}\cdots \GL_{n_s}$, where all $n_i$ are
even (since $u_1$ and $u_2$ are not conjugate under $\SO_n$).  This
implies that $L_1$ is not conjugate to $L_2$ under $\SO_n$.  
\end{proof}

\begin{proof}[Proof of Theorems \ref{theorem:map_with_regular_blocks},
  \ref{theorem:surjective_map_with_distinguished_parabolics}, and
  \ref{theorem:inverse_of_map_with_distinguished_parabolics}]
  
  There is essentially nothing to show for the case $G = \GL_n$, so we
  assume now that $G$ equals one of $\SO_n$ or $\Sp_n$.  
  
  Let $C$ be a unipotent $G$-class, $u\in C$ a unipotent element, let
  $\GL_{n_1} \cdots \GL_{n_s}H_m$ be a minimal Levi subgroup
  containing $u$ with $H$ equal to $\SO_{2m}$, $\SO_{2m+1}$ or
  $\Sp_{2m}$ as appropriate.  We factor $u$ as $u = u_1 \cdots u_s
  u_0$ with $u_i \in \GL_{n_i}$ for $1\le i \le s$ and $u_0\in
  \Cl_m^\circ$.  Note that $L$ has the maximal number of factors of
  the form $\GL_{n_i}$ among elements of $\R(G)$ which
  contain $u$, as required for $\Phi_1$; that $u_i$ is regular in
  $\GL_{n_i}$, for $i\ge 1$; and that $u_0$ is distinguished in
  $\Cl_m^\circ$.  Let $\alpha$ be the Jordan blocks of $u_1\cdots u_s$
  in the natural module for $\GL_{n_1}\cdots \GL_{n_s}$ and let
  $\beta$ be the Jordan blocks of $u_0$ in the natural module for
  $\Cl_m^\circ$.  Then the parts $\alpha$ equals $(n_1$, \dots,
  $n_s)$, $\beta$ satisfies the properties described in
  Lemma~\ref{lemma_mapping_subgroups_to_classes} parts (ii) and (iii),
  and, if $G=\SO_n$, $p=2$ and $n$ is even, then $\beta$ has an even
  number of parts.  These observations about $\beta$, combined with
  Lemma~\ref{lemma_mapping_subgroups_to_classes}, (and with the
  remaining part of this proof) also establish the assertions in
  Theorem~\ref{theorem:map_with_regular_blocks} about the image of
  $\Phi_1$.

  Note that $u = u_1\cdots u_s$ is equivalent to $\beta =0$ and
  $u=u_0$ is equivalent to $\alpha =0$.  Applying
  Lemma~\ref{lemma_mapping_subgroups_to_classes} it suffices to
  construct $\Phi_1$ and $\Phi_2$ under the assumption that $\alpha=0$
  or $\beta=0$.  In other words, if $\Phi_1$ and $\Phi_2$ have been so
  constructed, then the class represented by $\Phi_i(u)$ equals the
  class represented by $\Phi_i(u_1\cdots u_s) \Phi_i(u_0)$.
\begin{details}
  Let $H$ be the group so defined.  Then $\Psi_i(H)$ has the right
  Jordan blocks, so it suffices to verify that these blocks have the
  right singularity.  When $p\ne 2$ there is nothing to show, so we
  assume now that $p=2$.  Let $x$ be a part.  If $x$ is odd then it is
  singular wherever it occurs, thus there is nothing to show, so we
  assume now that $x$ is even.  If $x$ has odd multiplicity then it
  will occur as a part of $u_0$, and it will be nonsingular with
  respect to $u_0$ and with respect to $u$.  If $x$ has even
  multiplicity and is singular, then all parts of $\lambda(u)$ equal
  to $x$ will occur with respect to $u_1\cdots u_s$ and will be
  singular.  If $x$ has even multiplicity and is nonsingular, then
  exactly two parts of $\lambda(u)$ which equal $x$ will be parts of
  $\lambda(u_0)$ and will thus be nonsingular.
\end{details}

If $\beta=0$ we define $\Phi_i(u) = \GL_{n_1}\cdots \GL_{n_s}$.
  
Suppose $\alpha =0$.  For $p\ne 2$ or for $G = \Sp_{2n}$ we define
$\Phi_1(u) = \Cl_{m_1}^\circ\cdots \Cl_{m_r}^\circ$.  For $p=2$ and
$G=\SO_n$ we define $\Phi_1(u) = \Cl_{m_1}\cdots \Cl_{m_r}$.

For $\Phi_2(u)$ we define $\beta = \beta\one \oplus \beta\two \oplus
\beta\three$ as in
Theorem~\ref{theorem:inverse_of_map_with_distinguished_parabolics}.
Then each $\beta^{(i)}$ satisfies the properties in
Table~\ref{table_jordan_blocks_dist_rich_classes} (we prove this for
the case $G= \SO_n$ and $p=2$ in
Proposition~\ref{decomposing_distinguished_classes_into_Richardson_classes}
below) and each part $x$ satisfies $\epsilon(x) = 1$.  Therefore one
may apply the formulas in
Table~\ref{table:levi_partitions_distinguished} to find a unique
parabolic subgroup $P_i$ of $\Cl_{|\beta^{(i)}|}^\circ$ such that the
Jordan blocks of the Richardson class of $P_i$ equals $\beta^{(i)}$
(recall that Table~\ref{table:levi_partitions_distinguished} defines
an injective map from $G$-classes of distinguished parabolics to the
partitions of $n$ described in
Table~\ref{table_jordan_blocks_dist_rich_classes}).
\end{proof}

We say a partition $\beta$ satisfies the \define{difference
  condition} if for all even $i$ such that $\beta_{i+1}\ge 1$ we have
$\beta_i - \beta_{i+1} \ge 3$.
If $p=2$ and $\beta$ equals the Jordan blocks of the Richardson class
of a distinguished parabolic subgroup of an orthogonal group, then
$\beta$ satisfies the difference condition (see
Table~\ref{table_jordan_blocks_dist_rich_classes}).

\begin{prop}
\label{decomposing_distinguished_classes_into_Richardson_classes}
Let $\beta$ be a partition with at most one part equal to $1$ and all
other parts even with multiplicity at most $2$.  If $\beta$ does not
have a part equal to $1$ then we assume that $\beta$ has an even
number of parts.  Apply
Definition~\textup{\ref{remark_defining_decomposition}} for the case
$p=2$ and $G=\SO_n$ to decompose $\beta$ as $\beta = \beta\one \oplus
\delta = \beta\one \oplus \beta\two \oplus \beta\three$.  The
following hold:
\begin{enumerate}
\renewcommand{\labelenumi}{(\roman{enumi})}
\item If $\beta$ satisfies the difference condition then $\beta =
  \beta\one$.
\item If $\beta^{(i)}$ does not contain $1$ then it has an even number of
  parts.
\item Each $\beta^{(i)}$ satisfies the difference condition.
\item $\beta\three =0$ if and only if it is possible for $\beta$ to be
  decomposed into two partitions each of which satisfies the
  difference condition
\end{enumerate}
\end{prop}

\begin{proof}
Recall that we apply $f$ to parts of $\beta$, and that when we do so
we refer only to $\beta_j$, but we keep track (implicitly) of $j$.
Thus, $f(\beta_j)$ depends not only upon the value $\beta_j$ but also
on $j$.

Part (i).  Suppose that $\beta$ satisfies the difference condition.
An inductive argument on $j$ shows that $f(\beta_j) =1$ for all $j$.
\begin{details}
For $j = 1$ we have that none of conditions (1)--(3) in the definition
of $f$ can apply.  Thus $f(\beta_1) =1$.  

Now suppose $j\ge 2$ and that $f(\beta_1) = \dots = f(\beta_{j-1})=1$
and let $\ell$ be as in
Definition~\ref{remark_defining_decomposition}.  Note that $\ell =
j-1$ and $i=0$.  If $\ell$ is odd then none of the conditions (1)--(3)
in the definition of $f$ can apply, whence $f(\beta_j)= 1$.  If $\ell$
is even, then since $\beta$ satisfies the difference condition and
$\ell = j-1$ we have $\beta_{j-1}-\beta_j \ge 3$, whence condition (1)
cannot apply.  Finally, note that $i=0$ so conditions (2) and (3)
cannot apply.  Thus $f(\beta_j)=1$.
\end{details}

In the remaining proof, we use the phrase ``$\beta_j$ is added to
$\beta\one$'' to mean $f(\beta_j)=1$.  We also use obvious variations
on this.

Part (ii).  Since this decomposition is defined by applying $f$
recursively, it suffices to prove this claim for $\beta\one$ and
$\delta$.  Condition (2) in the definition of $f$ guarantees that the
last part of $\beta$ will be added to $\beta\one$ or $\delta$ in such
a way that both have an even number of parts, or that both have an
even number of parts greater than $1$.  This proves part (ii).  

In the remaining proof, we will use the notation ``$\beta_j\mapsto
\beta\one_\ell$'' to mean that $f(\beta_j)=1$ and that $\beta_j$
becomes the $\ell^{\rm th}$ part of $\beta\one$.  (Another way to say
this is that $\ell$ parts from $\beta_1$, \dots, $\beta_j$ map to $1$,
or, equivalently, map to $\beta\one$.)
  
Claim (a).  Let $\ell$ be even.  Then $\beta_j \mapsto \beta_\ell\one$
if and only if $\beta_{j-1} \mapsto \beta_{\ell-1}\one$.  Proof: This
follows immediately from the definition of $f$.
\begin{details}
Suppose $\beta_{j-1}\mapsto \beta_{\ell-1}\one$.  Then $f(\beta_j)=1$
since none of the conditions (1)--(3) apply.  

We prove the other direction by contradiction.  Suppose $\beta_j
\mapsto \beta_\ell\one$ and $f(\beta_{j-1})=0$.  Then $\ell-1$, parts
from $\beta_1$, \dots, $\beta_{j-2}$ map to $1$.  Since $\ell-1$ is
odd the definition of $f$ gives $f(\beta_{j-1})=1$, a contradiction.
\end{details}

We say that $\beta$ has a \define{bad sequence} if there exists an
even number $i$ such that $\beta_i = \beta_{i+1} > \beta_{i+2} =
\beta_{i+3}$ with $\beta_{i+1} - \beta_{i+2} = 2$ and $\beta_{i+2}\ne
0$.  
  
Claim (b).  If $\delta$ violates the difference condition then $\beta$
has a bad sequence.  
Proof: Let $i$ be even with $\delta_{i+1} \ge 1$ and $\delta_i -
\delta_{i+1}\le 2$.  Let $\beta_j$ map to $\delta_i$, let $\ell$ and
$k$ be as in the definition of $f$ (i.e.\ $\ell$ is the number of
parts from $\beta_1$, \dots, $\beta_{j-1}$ which map to $1$; $k$ is
the last of these parts which map to $1$; note that $i$ in the
definition of $f$ has been replaced by $i-1$ in the present context).
Then we have $\beta_k \mapsto \beta_\ell\one$.  Since $f(\beta_j) = 0$
we have that $\ell$ is even, whence $j= \ell + i$ is even

Suppose $\delta_i = \delta_{i+1}$.  Then $\delta_i \ge 2$, 
\begin{details}
Because at most one part of $\beta$ equals $1$.  
\end{details}
$\beta_j = \beta_{j+1}$ and $\beta_{j+1}$ maps to $\delta_{i+1}$.
\begin{details}
because $\beta_j \mapsto \delta_i$, and \emph{some} part of $\beta$
 maps to $\delta_{i+1}$, and this has to be the next part of
 $\beta$. 
\end{details}
Since $i$ is even and $f(\beta_{j+1}) =0$, we have that
$\beta_\ell\one - \beta_{j+1} = 2$,
\begin{details}
  Since $i$ is even it has to be condition (1) in the definition of
  $f$ which makes $f(\beta_{j+1})=0$.  Therefore $\beta_\ell\one -
  \beta_{j+1} \le 2$.  Since $\delta_i\ge2$ both these numbers are
  even, so the difference is $0$ or $2$.  We can't have $0$ since the
  multiplicity of all parts is at most $2$ and we already have
  $\beta_j = \beta_{j+1}$.
\end{details}
whence $k$ equals $j-1$ or $j-2$ (i.e.  $\beta_{j-1}$ or $\beta_{j-2}$
maps to $\beta_\ell\one$).  
\begin{details}
Because $\beta_{j-3} - \beta_j \ge 4$.  
\end{details}

We summarize this information below as
follows.  Each part $\beta_a$ is sent to a part $\beta\one_b$ in
$\beta\one$ or a part $\delta_b$ in $\delta$.  We indicate this by
writing $\underset{\textstyle \delta_b}{\beta_a}$ or
\raisebox{\height}{$\underset{\textstyle\beta_a}{\beta\one_b}$}
respectively.  We write $a \overset 2>b$ to indicate that $a-b=2$.
With this notation we have:
\[
\overset{\textstyle\beta_\ell\one}{\beta_{j-1}} 
\overset 2>
\underset{\textstyle\delta_i}{\beta_j} 
=
\underset{\textstyle\delta_{i+1}}{\beta_{j+1}}\,,\text{ or}\quad 
\overset{\textstyle\beta_\ell\one}{\beta_{j-2}}
=
\underset{\textstyle\delta_{i-1}}{\beta_{j-1}}
\overset 2>
\underset{\textstyle\delta_i}{\beta_j}
=
\underset{\textstyle\delta_{i+1}}{\beta_{j+1}}\,.
\]
The second subsequence in $\beta$ is a bad sequence; whence in this
case we are done proving claim (b).  In the first subsequence we have,
by claim (a), that $\beta_{j-2} \mapsto \beta_{\ell-1}\one$, 
\begin{details}
 since $\ell$ is even
\end{details}
thus, an even number of parts from $\beta_1$, \ldots, $\beta_{j-3}$
map to $1$
\begin{details}
this number is $\ell -2$
\end{details}
and an odd number of these parts map to $0$ (recall that $\ell$ is
even and $j$ is even).  
\begin{details}
Thus, the total number of parts is $j-3$, which is odd.  Out of these
$\ell-2$, which is even, map to $1$, and the rest map to $0$. 
\end{details}
Since $f(\beta_{j-2})=1$ we have that $\beta_{j-2} = \beta_{j-1}$ (if
$\beta_{j-2}>\beta_{j-1}$ then condition (3) in the definition of
$f$ would have caused $f(\beta_{j-2})=0$).
\begin{details}
  Suppose $\beta_{j-2}>\beta_{j-1}$ and consider defining
  $f(\beta_{j-2})$.  At this point an even number, $\ell-2$, of parts
  have gone to $\beta\one$.  At this point an odd number, $i-1$, of
  parts have gone to $\delta$.  We have (i) $\beta_{j-1} - \beta_{j+1}
  \le 2$, (ii) $\beta_{j+1} \ne 0$, (iii) $\beta_{j-2} -
  \beta_{j+1}\ge 3$.  
\end{details}
This means that the first subsequence is $\beta_{j-2} = \beta_{j-1}
\overset {2}{>} \beta_j = \beta_{j+1}$ which is a bad sequence.

We suppose now that $\delta_i - \delta_{i+1}$ equals $1$ or $2$ and
show that this leads to a contradiction.  Case 1: $\beta_{j} =
\beta_{j+1}$.  It is easy to show that conditions (2) or (3) could not
have caused $f(\beta_j) = 0$,
\begin{details}
  Since $\beta_{j+1} \ge 2$ we cannot have condition (2) rejecting
  $\beta_{j}$.

Now we show that condition (3) leads to a contradiction.  If it was
condition (3) which caused $f(\beta_j)=0$ then we would have
$\beta_{j+1}-\beta_{j+3} \le 2$.  Note $\beta_j = \beta_{j+1}$, and
$\beta_{j+3} >0$, and recall that at most one part can equal $1$.
Combining these observations with the sequence
$$
\beta_j = \beta_{j+1} > \beta_{j+2} \ge  \beta_{j+3} >0
$$
shows that $\beta_j - \beta_{j+3} = 2$, but then condition (3) does
not apply.
\end{details}
whence we have that $\beta_\ell\one - \beta_j = 2$.  Thus
$f(\beta_{j+1}) = 0$ 
\begin{details}
since condition (1) applies again to $\beta_{j+1}$
\end{details}
and $\delta_{i} = \delta_{i+1}$, a contradiction.

Case 2: $\beta_{j} > \beta_{j+1}$ and $\beta_\ell\one > \beta_j$.
Then $\beta_\ell\one - \beta_{j+1}\ge 3$
\begin{details}
since we have
$$
\beta_\ell\one > \beta_j > \beta_{j+1}
$$
where only the second inequality can differ by $1$ and the other must
differ by $2$. 
\end{details}
whence, $f(\beta_{j+1}) = 1$, 
\begin{details}
Since $\beta_\ell\one -\beta_{j+1}\ge 3$ we have that condition (1)
cannot reject $\beta_{j+1}$.

Since $i$ is even we have that conditions (2) and (3) cannot reject
$\beta_{j+1}$.
\end{details}
whence, by claim (a), we have that $f(\beta_{j+2}) = 1$. 
\begin{details}
Since $\beta_{j+1}\mapsto \beta_{\ell+1}\one$ with $\ell+1$ odd.
\end{details}
This implies that $\delta_i - \delta_{i+1} > 2$ contrary to
assumption.
\begin{details}
Since the first element of $\beta$ which can map to $\delta_{i+1}$ is
now $\beta_{j+3}$ and we have 
$$
\beta_j > \beta_{j+1} \ge \beta_{j+2} \ge \beta_{j+3} >0
$$
\end{details}

Case 3: $\beta_j > \beta_{j+1}$ and $\beta_\ell\one = \beta_j$.  We
have $\beta_{j-1} \mapsto \beta_\ell\one$, 
\begin{details}
Somebody before $\beta_j$ maps to $\beta_{\ell}\one$, but because all
multiplicities are $\le 2$, at can't be anyone before $\beta_{j-1}$.  
\end{details}
whence, by claim (a), we have that $\beta_{j-2} \mapsto
\beta_{\ell-1}\one$.  
\begin{details}
Because $\ell$ is even.
\end{details}
This shows that $f$ maps an even number of the parts $\beta_1$,
\ldots, $\beta_{j-3}$ to $1$ and an odd number of these parts to $0$.
\begin{details}
Specifically, $\ell-2$ parts map to $1$ and $i-1$ parts map to $0$.  
\end{details}
Note that $\beta_{j-1}-\beta_{j+1}\le 2$,
\begin{details}
Since $\delta_i - \delta_{i+1}$ equals $1$ or $2$ we have that
$\beta_j - \beta_{j+1}$ equals $1$ or $2$
\end{details}
that $\beta_{j+1} \ne 0$,
\begin{details}
Since $\delta_{i+1}\ge 1$
\end{details}
and that $\beta_{j-2} - \beta_{j+1} \ge 3$.
\begin{details}
since $\beta_{j-2} - \beta_{j-1} \ge 2$ and $\beta_j-\beta_{j+1}
=\{1,2\}$.  
\end{details}
But this would imply that $f(\beta_{j-2})=0$,
\begin{details}
by condition (3)
\end{details}
a contradiction.

This finishes the proof of claim (b).  

Claim (c). The partition $\delta$ has no bad sequence.  Proof:
Suppose, for contradiction, that $i$ is even, $\delta_i = \delta_{i+1}
\overset {2}{>} \delta_{i+2} = \delta_{i+3}$ with $\delta_{i+3}\ne 0$.
Let $\beta_j\mapsto \delta_i$.  Then we also have $\beta_{j+1}\mapsto
\delta_{i+1}$, $\beta_{j+2}\mapsto \delta_{i+2}$ and $\beta_{j+3}
\mapsto \delta_{j+3}$.  
\begin{details}
\emph{some} parts of $\beta$ have to map to $\delta_{i+1}$,
$\delta_{i+2}$, $\delta_{i+3}$, and comparing multiplicities it has to
be the next three parts of $\beta$ after $\beta_j$.  
\end{details}
Then condition (1) in the definition of $f$ cannot apply to either of
$\beta_{j+2}$ and $\beta_{j+3}$.  But conditions (2) and (3) both
require that an odd number of parts have already been mapped to
$\delta$.  Thus, it is not possible for both
$\beta_{j+2}$ and $\beta_{j+3}$ to be affected by conditions (2) and
(3), a contradiction.\smallskip

Part (iii).  The definition of $f$ makes it clear that $\beta\one$ and
$\beta\two$ satisfy the difference condition.  By claim (c), $\delta$
has no bad sequence, whence, by claim (b), $\beta\three$ satisfies the
difference condition.

Part (iv).  \forwards follows from part (iii).  Conversely, by claim
(b), it suffices to show that if $\beta$ has a bad sequence then it
cannot be written as the sum of two partitions each of which satisfies
the difference condition.  Suppose that $\beta$ has a bad sequence
$\beta_j = \beta_{j+1} > \beta_{j+2} = \beta_{j+3}$.  Fix a
decomposition $(\beta_1,\ \dots, \beta_{j-1}) = \mu \oplus \nu$ where
$\mu$ has an odd number of parts.  At most one part of $\beta_j$,
$\beta_{j+1}$, $\beta_{j+2}$ and $\beta_{j+3}$ can be added to $\mu$
without violating the difference condition and at most two of these
parts can be added to $\nu$ without violating the difference
condition.
\end{proof}

%
%

\section{Unipotent classes in classical subgroups of exceptional groups}
\label{section_extending_bala_carter}
In this section we give two examples of
Theorem~\ref{theorem:inverse_of_map_with_distinguished_parabolics},
and then return to our main application, translating the labels for
unipotent classes in a classical Levi subgroup of an exceptional group
$G$ into the labels for unipotent classes in $G$.  

First we recall some Bala-Carter notation.  If $R$ is a type of root
system then we use ``$R$'' to label the regular class in the simple
group of type $R$.  If $R$ is of type $A_n$, $B_n$, $C_n$ or $D_n$
then we use $R(a_j)$ to denote the Richardson class of the
distinguished parabolic whose Levi subgroup has only simple factors of
rank $1$ and a total semisimple rank of $j$.
%

\begin{example}
\label{example:decomposing_with_f}
Here we illustrate the map $f$, and the decomposition $\beta =
\beta\one \oplus \beta\two \oplus \beta\three$ it gives rise to as
described in Definition~\ref{remark_defining_decomposition}.
Throughout we assume that $G = \SO_n$ and $p=2$.

Suppose that $\beta = (8,4,1)$ equals the Jordan block sizes of an
unipotent element in $\SO_{13}$.  We start with $f(8) = 1$.  Now
$\beta_k = 8$, $\ell = 1$ and $i=0$, therefore $f(4) = 1$.  Now
$\beta_k = 4$, $\ell = 2$ and $i=0$.  Since $4-1 > 2$ we have $f(1) =
1$.  Therefore $\beta = \beta\one$ and $\beta\two = \beta\three = 0$.
This agrees with the fact that the unipotent class with Jordan blocks
given by $(8,4,1)$ is parameterized by the Bala-Carter Theorem; it is
the class $B_6(a_2)$.
\begin{details}
The diagram for $B_6(a_2)$ is given by {\tt xxxoxo}.  Thus $n_1 = 1$,
$n_2 =1$, $n_3 = 2$, $n_4=1$ and $m=1$.  So $\Lambda = (1^2,2)\oplus
1$.  We have $\psi(1^3 = n_i) = (2,1^{2\cdot 3-2},0) = (2,1^4)$,
$\psi(2) = 2^2$ and $\psi(1=m) = 2$ whence $\psi(\Lambda) = (2^4,
1^4)$ and $\lambda  = (2^4,1^4)^* \oplus 1 = (8,4)\oplus 1$.  
\end{details}

Suppose that $\beta = (12, 12, 10, 8, 6, 6, 4, 2)$ and we have
$G=\SO_{60}$.  For convenience we keep track of the results by writing
an array: the middle row has the original partition $\beta$, the first
row contains the parts that $f$ maps $1$, and the third row contains
the parts that $f$ maps to $0$.  The calculation of $f$ proceeds
sequentially from left to right.
$$
\begin{array}{rcccccccc}
\text{map to 1}:     & 12 & 12 &    &   & 6 & 6 &   &  \\
\beta:               & 12 & 12 & 10 & 8 & 6 & 6 & 4 & 2\\
\text{map to 0}:  &    &    & 10 & 8 &   &   & 4 & 2
\end{array}
$$
For instance, to calculate $f(10)$ one has that $\ell =2$ and
condition (1) causes $f(10)=0$.  Similarly, condition (3) causes $f(8)
= 0$, and condition (2) causes $f(2)=0$.  The decomposition is $\beta
= (12, 12, 6, 6)\oplus (10, 8, 4, 2)$.  Each of these partitions
corresponds to the Richardson class of a distinguished parabolic
subgroup of a group of type $D$.  The usual Bala-Carter notation does
not apply to these parabolics, but they can be described by the Dynkin
diagrams below (where each {\tt x} represent a crossed off node):
\begin{alltt}
                      o                             o
xo xo xoo xooo xooo x              x xo xooo xooo x  
                      o                             o
\end{alltt}
\begin{details}
For instance, consider (12,12,6,6)^* = (4,4,4,4,4,4,2,2,2,2,2,2) =
(4^6,2^6).  Since the largest part here is 4, we use the formula for
$p=2$ and $m_0>0$.  The largest part is $4=2m_0$, whence $m_0=2$.
The multiplicity of $4$ is $6$ whence $6=2c(4)+2$ and $c(4) = 2$.
There are no $3$'s so we have $2c(3)-2=0$ and $c(3) = 1$.  There are
$6$ $2$'s so we have $2c(2)+2=6$ whence $c(2)=2$.  Finally, there are
no $1$'s so we have $2c(1)-2=0$ and $c(1)=1$.  Thus, the Levi factors
are 
$$
\GL_1\GL_2\GL_2\GL_3\GL_4\GL_4\SO_{4}
$$
The labelled dynkin diagram is 
\begin{alltt}
                      o
xo xo xoo xooo xooo x
                      o
\end{alltt}
But this distinguished parabolic is not described through the usual
Bala-Carter notation.  It would be something like $D(b_{12})$ where
``$12$'' denotes semisimple rank 12, and ``$b$'' indicates that some
of the simple factors have rank $\ge 3$.  But this is ambiguous
because $D(b_{12})$ could also mean 
\begin{alltt}
                      o
xo xoo xoo xoo xooo x 
                      o
\end{alltt}

Similarly, the partition $(10,8,4,2)$ leads to $(10,8,4,2)^* =
(4,4,3,3,2,2,2,2,1,1) = (4^2,3^2,2^4,1^2)$.  This means that the
largest part is $4$, which equals $2m_0$, whence $m_0=2$.  The
multiplicity of $4$ is $2$ so $2=2c(4)+2$ whence $c(4)=2$.  Now
$mult(3)=2=2c(3)-2$ whence $c(3)=2$.  Similarly $mult(2)=4=2c(2)+2$
whence $c(2)=1$.  Finally $mult(1)=2=2c(1)-2$ whence $c(1)=2$.  Thus,
the Levi factor looks like $\GL_1\GL_1\GL_2\GL_3\GL_3\SO_4$ and the
Dynkin diagram looks like this
\begin{alltt}
                 o
x xo xooo xooo x 
                 o
\end{alltt}
\end{details}

Finally, consider the partition $(6,4,4,2,2,1)$ in $\SO_{19}$.
Applying $f$ once gives $\beta\one = (6,4)$ and $\delta = (4,2,2,1)$.
Applying $f$ to $\delta$ gives $\beta\two = (4,2)$ and $\beta\three =
(2,1)$.  
\end{example}

In the following example and lemmas, we describe Levi subgroups using
notation which specifies only their Lie type.  Thus, a Levi subgroup
denoted by, for example, $B_3T_1$, has Lie type of $B_3$ and a
central, one dimensional torus $T_1$.

\begin{example}
\label{example:classes_in_D_8}
If $G=D_8$ and $p=2$, then there are five extra classes.  In
Table~\ref{table:extra_classes_D_8} we give the Jordan blocks, the
decomposition of $\beta$ and a Bala-Carter type label.
\begin{table}
\caption{Extra unipotent classes in $D_8$, $p=2$}
\label{table:extra_classes_D_8}
$$
\begin{array}{ccc}
\parbox{.33\textwidth}
     {\raggedright Jordan blocks  ($\epsilon(4) = \epsilon(2)=1$ in all cases)} 
& 
\beta = \beta\one\oplus \beta\two \oplus \beta\three
&
\parbox{.33\textwidth}
    {Bala-Carter type label determined by Theorems
  \ref{theorem:surjective_map_with_distinguished_parabolics} and
  \ref{theorem:inverse_of_map_with_distinguished_parabolics}}
\\
\hline
(8, 4, 2^2) & (8, 4, 2^2) = (8,4)\oplus(2,2)     & D_6(a_1)D_2 \\
(6, 4^2, 2) & (6, 4^2, 2) = (6,4) \oplus (4,2)    & D_5(a_1)D_3 \\
(6,4,2^2, 1^2) & (6,4,2^2) = (6,4)\oplus(2,2) & D_5(a_1)D_2 \\
(4^2,2^4) & (4^2, 2^2) = (4,4)\oplus (2,2)       & A_1D_4(a_1)D_2 \\
(4^2, 2^2, 1^4) & (4^2,2^2) = (4,4)\oplus (2,2) & D_4(a_1)D_2 \\
\end{array}
$$
\end{table}


We sketch how one can see that these classes are extra.  By
Lemma~\ref{lemma_mapping_subgroups_to_classes}, the class $(4^2,
2^2,1^4)$ is distinguished in the Levi subgroup $D_6T_2$.  Since it
does not satisfy the properties in
Table~\ref{table_jordan_blocks_dist_rich_classes} applied to $D_6$, it
cannot be in the image of the Bala-Carter map.

Similarly the class $D_5(a_1)D_2$ is distinguished in the Levi
subgroup $D_7T_1$, the class $A_1D_4(a_1)D_2$ is distinguished in the
Levi subgroup $A_1D_6T_1$, and the classes $D_6(a_1)D_2$ and
$D_5(a_1)D_3$ are distinguished in $D_8$.
\end{example}

We turn now to the exceptional groups $E_7$, $E_8$ and $F_4$ and use
Lawther~\cite{lawther2} (who draws on the work of \cite{mizuno1},
\cite{mizuno2}, \cite{shinoda}, \cite{shoji}) for the number of
unipotent classes, the number of extra unipotent classes, and their
representatives.  We note that $E_6$ has no extra classes in any
characteristic (see \cite{duckworth}).
%
%

\begin{lemma}
  Let $p=2$, $X=B_3T_1$, and $G=F_4$.  Then $X$ has two extra
  unipotent classes.  One of these is distinguished in a $B_2T_2$ Levi
  subgroup, and we denote this class by $D_2$.  The other is
  distinguished in $X$ and we denote this class by $D_3$.  The classes
  $D_2$ and $D_3$ equal the $F_4$-classes Lawther denotes by $\tilde
  A_1^{(2)}$ and $B_2^{(2)}$, respectively.
\end{lemma}


\begin{proof}
We proceed as in Examples~\ref{example:decomposing_with_f} and
\ref{example:classes_in_D_8} to find the extra classes and decompose
their partitions.    They are $(2^2,1^3) = (2^2)\oplus (1^3)$ and
$(4,2,1) = (4,2)\oplus (1)$.  The class $(2,2)$ is regular in $\SO_4$,
whence we label it as $D_2$.  The class $(4,2)$ is regular in $\SO_6$,
whence we label it as $D_3$.  
\begin{details}
As pointed out after
Theorem~\ref{theorem:inverse_of_map_with_distinguished_parabolics}, a
class is extra if and only if $\beta\ne \beta\one$ (in the notation of
the theorem).  It is a simple matter to list all the possible Jordan
blocks of unipotent classes in $B_3$, as well as the possible values
of $\epsilon$ for those classes.  It is also simple to list the Jordan
blocks and $\epsilon$ values for those classes produced by the
Bala-Carter theorem.  Comparing the two lists, one finds two extra
classes, with Jordan blocks $(2^2,1^3)$, $\epsilon(2) = 1$, and
$(4,2,1)$.

Consider first the class $\beta = (2^2,1^3)$, $\epsilon(2)=1$.  We
indicate now how to decompose the partition as described before
Theorem~\ref{theorem:inverse_of_map_with_distinguished_parabolics}.
We start with $f(\beta_1)=1$.  This makes $\ell = 1$ and so we have
$f(\beta_2)=1$.  Now $\ell = 2$, $i=0$, $\beta_k=2$ and
$\beta_k-\beta_3 \le 2$.  Thus, $f(\beta_3)=0$.  Now $\ell = 2$,
$i=1$, $\beta_k=2$ and $\beta_k-\beta_4 \le 2$.  Thus,
$f(\beta_4)=0$.  Similarly, $f(\beta_5)=0$.  Thus, the decomposition
is $\beta\one = (2,2)$ and $\delta = (1,1,1)$.  The class $(2,2)$ is
the regular class in $\SO_4$, which we label as $D_2$.  

Similarly, the class $(4,2,1)$ decomposes as $(4,2)$ and $(1)$.  The
class $(4,2)$ is the regular class in $\SO_6$, which we label as
$D_3$.  
\end{details}

The classes $D_2$ and $D_3$ are distinguished in the $B_2$ and $B_3$
Levi subgroups of $F_4$ and are not in the image of the Bala-Carter
map.  The same comment applies to the classes $\tilde A_1^{(2)}$ and
$B_2^{(2)}$.  Applying
Corollary~\ref{bijection_reduction_to_dist_classes} we conclude that
the $D_2$ and $\tilde A_1^{(2)}$ classes are equal, as are the $D_3$
and $B_2^{(2)}$ classes.
\end{proof}

\begin{lemma}
  Let $p=2$, $X=D_6T_1$, and $G=E_7$.  Then $X$ has one extra
  unipotent class.  It is distinguished in $X$ and denoted by
  $D_4(a_1)D_2$.  Using Lawther's notation this is the class $A_3 +
  A_2^{(2)}$ in $E_7$.
\end{lemma}

\begin{proof}
  If $G=E_7$ and $p=2$ then there is one extra class.  The Levi
  subgroup $D_6T_1$ also has one extra class (see
  Example~\ref{example:classes_in_D_8}), which is distinguished in
  $D_6T_1$.  By Corollary~\ref{bijection_reduction_to_dist_classes}
  these extra classes are the same class, whence we label it as
  $D_4(a_1)D_2$.  Lawther~\cite{lawther2} denotes this class by $A_3 +
  A_2^{(2)}$.
\end{proof}

\begin{lemma}
  Let $p=2$, $X=D_7T_1$ and $G=E_8$.  Then $X$ has $2$ extra classes.
  One of these is the class $D_4(a_1)D_2$ described in the previous
  lemma.  The other is distinguished in $X$ and denoted by
  $D_5(a_1)D_2$.  Using Lawther's notation these are the classes
  $A_3+A_2^{(2)}$ and $D_4+A_2^{(2)}$ respectively.
\end{lemma}

\begin{proof}
  Example~\ref{example:classes_in_D_8} shows that $X$ has two extra
  classes and the previous lemma shows that $D_4(a_1)D_2$ and
  $A_3+A_2^{(2)}$ are the same class.  Thus, it remains to show that
  $D_5(a_1)D_2$ and $D_4+A_2^{(2)}$ are the same class.  
  
  Here is one way to verify this.  In the natural module for $D_7$ the
  class $D_5(a_1)D_2$ has Jordan blocks given by $(6,4,2,2)$.
  Decomposing these blocks as $(6,2)\oplus (4,2)$ we see that this
  class can be represented by a regular unipotent element in a
  $\SO_{8}\SO_6$ subgroup.  By \cite{steinberg_regular_elements} such
  a regular element can be represented as a product of elements in
  root groups with the roots forming a basis for a $D_4 + D_3$
  subsystem.  The roots given in \cite{mizuno2} for the
  $D_4+A_2^{(2)}$ class form the basis for a $D_4 + D_3$ root system.
\begin{details}
Mizuno numbers the roots as 40, 39, 37, 104, 41, 38 and 108.  They
are given by 
0122100
  1

1111110
  1

1221000
  1

0001111
  0

0121110
  1

1121100
  1

1111111
  0

respectively.

Let's label these roots as M40, M39. M37, M104, M41, M38 and M108.
Then we have
                M41
               /
    M37----M104
               \
                M38

and 
    M39---M40---M108
\end{details}
In $E_8$ all $D_4+D_3$ root systems are conjugate, so the root groups
used to represent the class we have called $D_5(a_1)D_2$ can be
conjugated to the roots given by Mizuno.
\end{proof}

\end{document}